\title{\bf{Some algebraic surfaces with canonical map of degree 10, 12, 14}}
\author{
	NGUYEN BIN\\
}
\date{\today}

\newcommand{\CurrentAddresses}{{
		\bigskip
		\footnotesize
		\text{Current address:}\par\nopagebreak
		\text{Department of Mathematics and Statistics,}\par\nopagebreak	
		\text{Quy Nhon University,}\par\nopagebreak	
		\text{Vietnam.}\par\nopagebreak		
		\textit{E-mail address}: \texttt{nguyenbin@qnu.edu.vn}
		
}}

\newcommand{\Addresses}{{
		\bigskip
		\footnotesize
		    \text{National Center for Theoretical Sciences,}\par\nopagebreak	
			\text{National Taiwan University,}\par\nopagebreak	
			\text{Taiwan.}\par\nopagebreak		
		\textit{E-mail address}: \texttt{nguyenbin@ncts.ntu.edu.tw}
				
	}}

\newcommand{\SecondAddresses}{{
		\bigskip
		\footnotesize
		\text{Center for Mathematical Analysis, Geometry and Dynamical Systems}\par\nopagebreak
		\text{Departamento de Matem\'{a}tica}\par\nopagebreak
		\text{Instituto Superior T\'{e}cnico}\par\nopagebreak
		\text{Universidade de Lisboa}\par\nopagebreak
		\text{Portugal.}\par\nopagebreak
		\textit{E-mail address}: \texttt{nguyenbin@tecnico.ulisboa.pt}
		
}}

\newcommand\blfootnote[1]{%
	\begingroup
	\renewcommand\thefootnote{}\footnote{#1}%
	\addtocounter{footnote}{-1}%
	\endgroup
}

\date{\today}

\documentclass[10pt]{article}
\usepackage{amsmath, amsthm, amssymb, mathrsfs, amscd, amsthm, amscd,amsfonts,calligra,mathrsfs,adjustbox,lipsum}

\usepackage{indentfirst,url,xypic}
\usepackage[all]{xy}
\usepackage{url}
\usepackage{tikz}
\usepackage[colorlinks,plainpages]{hyperref}
\hypersetup{
	colorlinks=true,
	linkcolor=blue,
	filecolor=magenta,      
	urlcolor=cyan,
}

\DeclareMathOperator{\degree}{deg}
\DeclareMathOperator{\image}{im}

\newtheorem{Theorem}{Theorem}
\newtheorem{Proposition}{Proposition }

\newtheorem{Notation}[Proposition]{Notation}

\theoremstyle{remark}
\newtheorem{Remark}{Remark}

\newcommand{\MSC}{\textbf{Mathematics Subject Classification (2010):}}

\newcommand{\Key}{\textbf{Key words:}}

\begin{document}
\maketitle
\begin{abstract}
	Surfaces of general type with canonical map of degree $ d $ bigger than $ 8 $ have bounded geometric genus and irregularity. In particular the irregularity is at most $ 2 $ if $ d\geq 10 $. In the present paper, the existence of surfaces with $ d=10 $ and all possible irregularities, surfaces with $ d = 12 $ and irregularity $ 1 $ and $ 2 $, and surfaces with $ d = 14 $ and irregularity $ 0 $ and $ 1 $ is proven, by constructing these surfaces as $ \mathbb{Z}_2^3 $-covers of certain rational surfaces. These results together with the construction by C. Rito of a surface with $ d=12 $ and irregularity $ 0 $ show that all the possibilities for the irregularity in the cases $ d=10 $, $ d=12 $ can occur, whilst the existence of a surface with $ d=14 $ and irregularity $ 2 $ is still an open problem.
\end{abstract}

\blfootnote{\MSC{ 14J29}.}
\blfootnote{\Key{ Surfaces of general type, Canonical maps, Abelian covers.}}

\section{Introduction} 
 The problem of constructing examples of surfaces of general type with canonical map of high degree has been studied by many authors in the last decades. Let $ X $ be a minimal smooth complex surface of general type and denote by $ \xymatrix{\varphi_{\left| K_X\right| }:X \ar@{.>}[r] & \mathbb{P}^{p_g\left( X\right)-1}} $ the canonical map of $ X $, where $ K_{X} $ is the canonical divisor of $ X $ and $ p_g\left( X\right) = \dim H^0\left( X, K_{X}\right)  $ is the geometric genus. In 1979, A. Beauville proved that the degree $ d $ of the canonical map is at most $ 9 $ if the surface has holomorphic Euler-Poincar\'{e} characteristic bigger than $ 30 $ {\cite[ \rm Proposition 4.1]{MR553705}}. Later in 1986, G. Xiao showed that the degree of the canonical map is at most $ 8 $ if the surface has geometric genus bigger than $ 132 $ {\cite[ \rm Theorem 3]{MR842626}}. Only few surfaces with $ d $ greater than $ 8 $ have been known so far, such as: S. L. Tan's example \cite{MR1141782} with $ d = 9 $, U. Persson's example \cite{MR527234} with $ d = 16 $. In the last decade, some surfaces with $ d = 12, 16, 20, 24,27,32,36 $ were constructed by C. Rito \cite{MR3391024}, \cite{MR3619737}, \cite{MR3663791}, \cite{MR4372410}, C. Gleissner, R. Pignatelli and C. Rito \cite{2018arXiv180711854G}, Ching-Jui Lai and Sai-Kee Yeung \cite{MR4298919}, \cite{2015arXiv151006622Y} and the author \cite{MR4008073}, \cite{MR4334862}. There are recent preprints \cite{2022arXiv220702969F}, \cite{2022arXiv220906057F} of F. Fallucca and C. Gleissner constructing surfaces with $ d= 10, 11, 12, 13, 14, 15, 18 $.\\
 
 When the canonical map has degree $d = 10, 12, 14 $, by the Bogomolov-Miyaoka-Yau inequality the irregularity $ q $ is at most $ 2 $. Indeed, we have 
 \begin{align*}
 	d\left(p_g -2 \right)   \le d\degree\left( \image\left( \varphi_{\left| K_X\right| }\right) \right) \le K_X^2 \le 9\chi\left( \mathcal{O}_X\right)  = 9p_g -9q +9.
 \end{align*}
 \noindent
 This implies that $ q \le 2 $ since $ p_g \ge 3 $. In this paper, by using $ \mathbb{Z}_2^3 $-covers of del Pezzo surfaces of degree $ 5 $, we construct three surfaces with $ d = 10 $ and with all possible irregularities $ q = 0,1,2 $. By using $ \mathbb{Z}_2^3 $-covers of del Pezzo surfaces of degree $ 6 $, we construct two surfaces with $ d = 12 $ and with $ q = 1,2 $. These two examples together with C. Rito's with $ d = 12 $ and with $ q = 0 $ show the existence of surfaces with $ d = 12 $ for every possibility of the irregularity. By using $ \mathbb{Z}_2^3 $-covers of del Pezzo surfaces of degree $ 7 $, we construct two surfaces with $ d = 14 $ and with $ q = 0,1 $. We do not have any example with $ d = 14 $ and with $ q = 2 $. The following theorem is the result of this paper:
 \begin{Theorem}\label{the theorem with d = 10 12 14}
 	There exist minimal surfaces of general type $ X $ satisfying the following
 	$$
 	\begin{tabular}{|c |c| c| c| c|}
 		\hline
 		$ d $&$ q\left( X\right) $ &$ p_g\left( X\right) $ &$ K_X^2 $& $\left| K_X \right| $ \\
 		\hline
 		$ 10 $&$ 0 $&$ 3 $&$ 10 $& is base point free\\
 		\hline
 		$ 10 $&$ 1 $&$ 3 $&$ 10 $& is base point free\\
 		\hline
 		$ 10 $&$ 2 $&$ 3 $&$ 12 $& has a non-trivial fixed part\\
 		\hline
 		$ 12 $&$ 0 $&$ 3 $&$ 12 $& is base point free\\
 		\hline
 		$ 12 $&$ 1 $&$ 3 $&$ 12 $& is base point free\\
 		\hline
 		$ 12 $&$ 2 $&$ 3 $&$ 12 $& is base point free\\
 		\hline
 		$ 14 $&$ 0 $&$ 3 $&$ 14 $& is base point free\\
 		\hline
 		$ 14 $&$ 1 $&$ 3 $&$ 14 $& is base point free\\
 		\hline
 	\end{tabular} 
 	$$
 	\noindent 
 	Furthermore, the canonical map corresponding to each surface in the above table is a morphism.
 \end{Theorem}   
  The idea of our constructions is the following: we first construct the surfaces with $ d = 14 $ and with $ q = 0, 1 $ by taking specific $ \mathbb{Z}_{2}^3 $-covers $ X $ of del Pezzo surface $ Y_2 $ of degree $ 7 $ (see Notation \ref{Notation of del Pezzo surface of degree 7}). In particular, we take the covers to be in such a way that the following diagram commutes
 $$
 \xymatrix{X \ar[0,3]^{\mathbb{Z}_2^3}_g \ar[1,1]^{\mathbb{Z}_2} \ar@{.>}[2,1]_{\varphi_{\left| K_X \right| }}&&& Y_2\\
 	&Z \ar[-1,2]_{\mathbb{Z}_2^2} \ar[1,0]^{\varphi_{\left| K_{Z} \right| }}&&\\
 	&\mathbb{P}^{2}&&} 	
 $$
 \noindent
 where the intermediate surface $ Z $, which is a $ \mathbb{Z}_{2} $-quotient of $X$, is a surface of general type whose only singularities are of type $ A_1 $. Moreover, we require that: 
 \begin{enumerate}
 	\item the canonical map of $ Z $ is of degree $ 7 $,
 	\item the canonical map of $ X $ factors through the quotient map $ \xymatrix{X \ar[r] & Z} $.	
 \end{enumerate}
 \noindent
 More precisely, we take the $ \mathbb{Z}_{2}^3 $-covers such that the $ \mathbb{Z}_{2} $-quotient $ Z:=X/\left( 0,0,1\right)  $ is the bidouble cover of $ Y_2 $ with the building data $ \left\lbrace D_1, D_2, D_3, L_1, L_2, L_3\right\rbrace  $ determined as follows:
 \begin{align*}
 	D_1:= &D_{010} + D_{011} \equiv -K_{Y_2}\\ 
 	D_2:= &D_{100} + D_{101} \equiv -K_{Y_2}\\
 	D_3:= &D_{110} + D_{111} \equiv -K_{Y_2}\\
 	L_1:= &L_{100} \equiv -K_{Y_2}\\
 	L_2:= &L_{010} \equiv -K_{Y_2}\\
 	L_3:= &L_{110} \equiv -K_{Y_2}.
 \end{align*}
 \noindent
 The branch divisors $ D_{\sigma} $ are chosen so that the intermediate surface $ Z $ has only singularities of type $ A_1 $ and that the canonical linear system $ \left| K_Z\right|  $ is base point free. So the canonical map of $ Z $ is of degree $ 7 $.\\
 
 \noindent
 Furthermore, we take the $ \mathbb{Z}_{2}^3 $-covers with $ h^{0}\left( L_{\chi} + K_{Y}\right) = 0  $ for all $ \chi \in \left\lbrace  \chi_{001}, \chi_{101}, \chi_{011}, \chi_{111}  \right\rbrace  $ so that the canonical map of $ X $ factors through the quotient map $ \xymatrix{X \ar[r] & Z} $. In fact, by Proposition \ref{invariants of Z_2^n cover}, we have the following decomposition:     
 \begin{align*}
 	H^{0}\left( X, K_X\right) = H^{0}\left( Y_2, K_{Y_2}\right) \oplus \bigoplus_{\chi \ne  \chi_{000}}{H^{0}\left( Y_2, K_{Y_2} +L_{\chi}\right)}. 
 \end{align*}
 
 \noindent
 We consider the subgroup $ \Gamma:= \left\langle \left( 0,0,1\right) \right\rangle  $ of $ \mathbb{Z}_2^3 $. Let $ \Gamma^\perp $ denote the kernel of the restriction map $ \xymatrix{\left( \mathbb{Z}_2^3\right)^{*} \ar[r]&\Gamma^{*}} $, where $ \Gamma^{*} $ is the character group of $ \Gamma $. We have $ \Gamma^{\perp} = \left\langle \chi_{100},\chi_{010},\chi_{110} \right\rangle  $. The subgroup $ \Gamma $ acts trivially on $ H^{0}\left( X, K_X\right) $ since $ h^0\left( L_{ \chi} + K_{Y_2} \right)  = 0 $ for all $ \chi \notin \Gamma^\perp $. So the canonical map $ \varphi_{\left| K_X \right| } $ is the essential composition of the quotient map $ \xymatrix{X \ar[r]& Z:= X/\Gamma } $ with the canonical map $\varphi_{\left| K_{Z} \right| }$ of $ Z $ (see e.g. \cite[\rm Example 2.1]{MR1103913}, \cite[\rm Section 3.2]{MR4203293}, \cite[\rm Section 3.1.2]{MR4394091}).\\   
 
 The remaining surfaces described in Theorem \ref{the theorem with d = 10 12 14} are obtained as degenerations of the two surfaces mentioned above. We modify these $ \mathbb{Z}_{2}^3$-covers by imposing triple and quadruple points to the branch loci. Subsequently, we resolve the resulting singularities. As part of our degeneration process, we create a $ T $-singularity, specifically a singular point of type $ \frac{1}{4}\left( 1,1\right) $. These singular surfaces are examples of $ \mathbb{Q} $- Gorenstein smoothing, where the singular surfaces stay in the boundary of the KSBA moduli space of surfaces of general type with fixed topological invariants $ \chi\left( \mathcal{O} \right) $ and $ K^2 $.

	\section{Notation and $ \mathbb{Z}_{2}^3$-coverings}
	\subsection{Notation and conventions}
	Throughout this paper all surfaces are projective algebraic over the complex numbers. Linear equivalence of divisors is denoted by $ \equiv $. A character $ \chi $ of the group $ \mathbb{Z}_{2}^3$ is a homomorphism from $ \mathbb{Z}_{2}^3$ to $ \mathbb{C}^{*} $, the multiplicative group of the non-zero complex numbers. We also use the following notation of del Pezzo surfaces:
	
	\begin{Notation}\label{Notation of del Pezzo surface of degree 7} 
	We denote by $ Y_2 $ the blow-up of $ \mathbb{P}^2$ at two distinct points $ P_1 $, $ P_2 $, by $ l $ the pull-back of a general line in $ \mathbb{P}^2$, by $ e_1 $, $ e_2 $ the exceptional divisors corresponding to $ P_1 $, $ P_2 $, respectively, by $ f_1 $, $ f_2$ the strict transforms of a general line through $ P_1 $, $ P_2 $, respectively and by $ h_{12} $ the strict transform of the line $ P_1  P_2 $. The anti-canonical class $ -K_{Y_2} \equiv f_1 + f_2 + l $ is very ample and the linear system $ \left| -K_{Y_2} \right|  $ embeds $ Y_2 $ as a smooth del Pezzo surface of degree $ 7 $ in $ \mathbb{P}^7 $.	
	\end{Notation}	   	     
	\begin{Notation}\label{Notation of del Pezzo surface of degree 6} 
		We denote by $ Y_3 $ the blow-up of $ \mathbb{P}^2$ at three distinct non-collinear points $ P_1 $, $ P_2 $, $ P_3 $. Let us denote by $ l $ the pull-back of a general line in $ \mathbb{P}^2$, by $ e_1 $, $ e_2 $, $ e_3 $ the exceptional divisors corresponding to $ P_1 $, $ P_2 $, $ P_3 $, respectively, by $ f_1 $, $ f_2$, $ f_3$ the strict transforms of a general line through $ P_1 $, $ P_2 $, $ P_3 $, respectively and by $ h_{12} $, $ h_{23} $, $ h_{31} $ the strict transforms of the lines $ P_1  P_2 $, $ P_2  P_3 $, $ P_3  P_1 $, respectively. The anti-canonical class $ -K_{Y_3} \equiv f_1 + f_2 + f_3 $ is very ample and the linear system $ \left| -K_{Y_3} \right|  $ embeds $ Y_3 $ as a smooth del Pezzo surface of degree $ 6 $ in $ \mathbb{P}^6 $.
	\end{Notation}
	
	\begin{Notation}\label{Notation of del Pezzo surface of degree 5}
	We denote by $ Y_4 $ the blow-up of $ \mathbb{P}^2$ at four points in general position $ P_1, P_2, P_3, P_4 $. Let us denote by $ l $ the pull-back of a general line in $ \mathbb{P}^2$, by $ e_1 $, $ e_2 $, $ e_3 $, $ e_4 $ the exceptional divisors corresponding to $ P_1 $, $ P_2 $, $ P_3 $, $ P_4 $, respectively, by $ f_1 $, $ f_2$, $ f_3$, $ f_4$ the strict transforms of a general line through $ P_1 $, $ P_2 $, $ P_3 $, $ P_4 $, respectively and by $ h_{ij} $ the strict transforms of the line $ P_i  P_j $, for all $ i, j \in \left\lbrace 1,2,3,4\right\rbrace  $, respectively. The anti-canonical class 
	\begin{align*}
	-K_{Y_4} &\equiv f_1 + f_2 + f_3 - e_4 \\
	&\equiv f_1 + f_2 + f_4 - e_3 \\
	&\equiv f_1 + f_3 + f_4 - e_2 \\
	&\equiv f_2 + f_3 + f_4 - e_1  
	\end{align*}
	\noindent
	is very ample and the linear system $ \left| -K_{Y_4} \right|  $ embeds $ Y_4 $ as a smooth del Pezzo surface of degree $ 5 $ in $ \mathbb{P}^5 $.    	   		
	\end{Notation}

	\subsection{$ \mathbb{Z}_{2}^3 $-coverings}	  
	An abelian cover of a smooth surface $ Y $ with finite abelian group $ G $ is a finite map $ \xymatrix{f: X \ar[r] & Y} $ together with a faithful action of $ G $ on $ X $ such that $ f $ exhibits $ Y $ as the quotient of $ X $ via $ G $.\\
	
	The construction of abelian covers was studied by R. Pardini in \cite{MR1103912}. For detail about the building data of abelian covers and their notations, we refer the reader to Section 1 and Section 2 of R. Pardini's work. For the sake of completeness, we recall some facts on $ \mathbb{Z}_{2}^3 $-covers, in a form which is convenient for our later constructions.\\

	    We will denote by  $ \chi_{j_1j_2j_3} $ the character of $ \mathbb{Z}_{2}^3 $ defined by
	    \begin{align*}
	    \chi_{j_1j_2j_3}\left( a_1,a_2,a_3\right): =  e^{\left( \pi a_1j_1\right)\sqrt{-1}}e^{\left( \pi a_2j_2\right) \sqrt{-1}}e^{\left( \pi a_3j_3\right) \sqrt{-1}}
	    \end{align*}
	    for all $ j_1,j_2,j_3,a_1,a_2,a_3\in \mathbb{Z}_2 $. A $ \mathbb{Z}_{2}^3 $-cover $ \xymatrix{X \ar[r] & Y} $ can be determined a collection of non-trivial divisors $ L_{\chi} $ labelled by characters of $ \mathbb{Z}_{2}^3 $ and effective divisors $ D_{\sigma} $ labelled elements of $ \mathbb{Z}_{2}^3 $ of the surface $ Y $. More precisely, from \cite[\rm Theorem 2.1]{MR1103912} we can define $ \mathbb{Z}_{2}^3 $-covers as follows:
	    \begin{Proposition} \label{Construction of cover of degree 8}
	    	Given $ Y $ a smooth projective surface with no non-trivial $ 2 $-torsion, let $ L_{\chi} $ be divisors of $ Y $ such that $ L_{\chi} \not\equiv \mathcal{O}_Y $ for all non-trivial characters $ \chi $ of $ \mathbb{Z}_{2}^3  $ and let $ D_{\sigma} $ be effective divisors of  $ Y $ for all $ \sigma \in \mathbb{Z}_{2}^3 \setminus \left\lbrace \left(0,0,0 \right)  \right\rbrace  $ such that the total branch divisor $ B:=\sum\limits_{\sigma \ne 0}{D_{\sigma}} $ is reduced. Then $ \left\lbrace L_{\chi}, D_{\sigma} \right\rbrace_{\chi,\sigma}$ is the building data of a $ \mathbb{Z}_{2}^3$-cover $ \xymatrix{f:X \ar[r]& Y} $ if and only if
	    	$$
	    	\begin{adjustbox}{max width=\textwidth}
	    	\begin{tabular}{l l l l l l l l }
	    	$ 2L_{100} $&$ \equiv $&$  $&$ $&$ D_{100} $&$ +D_{101 } $&$ +D_{110} $&$ +D_{111} $ \\
	    	$ 2L_{010} $&$ \equiv $&$ D_{010} $&$ +D_{011} $&$  $&$ $&$ +D_{110} $&$ +D_{111} $ \\
	    	$ 2L_{001} $&$ \equiv D_{001 } $&$ $&$ +D_{011} $&$$&$ +D_{101} $&$  $&$ +D_{111 } $ \\
	    	$ 2L_{110} $&$ \equiv $&$ D_{010 } $&$ +D_{011} $&$ +D_{100} $&$ +D_{101} $&$  $&$  $\\
	    	$ 2L_{101} $&$ \equiv D_{001} $&$ $&$ +D_{011} $&$ +D_{100} $&$  $&$ +D_{110} $&$  $ \\
	    	$ 2L_{011} $&$ \equiv D_{001} $&$ +D_{010} $&$  $&$  $&$ +D_{101} $&$ +D_{110} $&$  $ \\
	    	$ 2L_{111} $&$ \equiv D_{001} $&$ +D_{010} $&$  $&$ +D_{100} $&$ $&$ $&$ +D_{111 } $.
	    \end{tabular}
	\end{adjustbox}
	$$	
	\end{Proposition}
	
	\noindent
	The condition $ L_{\chi} \not\equiv \mathcal{O}_Y $ for all non-trivial characters $ \chi $ assures that the surface $ X $ is irreducible.\\
	
	By \cite[\rm Theorem 3.1]{MR1103912} if each branch component $D_\sigma$ is smooth and the total branch locus $B $ is a simple normal crossings divisor, the surface $X$ is smooth. \\
	
	Also from \cite[\rm Lemma 4.2, Proposition 4.2]{MR1103912} we have:
	\begin{Proposition}\label{invariants of Z_2^n cover}
	If $ Y $ is a smooth surface and $ \xymatrix{f: X \ar[r]& Y} $ is a smooth $  \mathbb{Z}_{2}^3$-cover with the building data $ \left\lbrace L_{\chi}, D_{\sigma} \right\rbrace_{\chi,\sigma}$, the surface $ X $ satisfies the following:
	\begin{align*}
		2K_X & \equiv f^*\left( 2K_Y + \sum\limits_{\sigma \ne 0} {D_{\sigma} } \right); \\
		f_{*}\mathcal{O}_X &= \mathcal{O}_Y \oplus \bigoplus\limits_{\chi \ne \chi_{000}  }L_{\chi}^{-1}.
	\end{align*}	
	\noindent
	This implies that	
	\begin{align*}
		H^{0}\left( X, K_X\right) &= H^{0}\left( Y, K_{Y}\right) \oplus \bigoplus_{\chi \ne  \chi_{000}}{H^{0}\left( Y, K_{Y} +L_{\chi}\right)};\\
		K^2_X &= 2\left( 2K_Y + \sum\limits_{\sigma \ne 0} {D_{\sigma} } \right)^2; \\
		p_g\left( X \right) &=p_g\left( Y \right) +\sum\limits_{\chi \ne  \chi_{000}  }{h^0\left( L_{\chi} + K_Y \right)}; \\
		\chi\left( \mathcal{O}_X \right) &= 8\chi\left( \mathcal{O}_Y \right)  +\sum\limits_{\chi \ne \chi_{000}  }{\frac{1}{2}L_{\chi}\left( L_{\chi}+K_Y\right)}. 
	\end{align*}
	\end{Proposition}
	\noindent
	We notice that most of Proposition \ref{invariants of Z_2^n cover} does not require the smoothness assumption on X. For a more general statement, we direct the reader to the work of Bauer and Pignatelli \cite[\rm Section 2]{MR4278662}.\\
	
	From \cite[\rm Proposition 4.1,c 2.1]{MR1103912}, the generators of the canonical linear system $ \left|  K_X \right|  $ are obtained as follows:
	\begin{Proposition}\label{The generators of canonical system of Z_2^n cover}
	If $ Y $ is a smooth surface and $ \xymatrix{f: X \ar[r]& Y} $ is a smooth $  \mathbb{Z}_{2}^3 $-cover with building data $ \left\lbrace L_{\chi}, D_{\sigma} \right\rbrace_{\chi,\sigma}$, the canonical linear system $ \left|  K_X \right|  $ is generated by
	\begin{align*}
	f^{*}\left(  \left| K_Y + L_{\chi}\right| \right) +\sum\limits_{\chi\left( \sigma\right)=1 }{R_{\sigma}}    , \hskip 5pt \forall \chi \in J 
	\end{align*}
	\noindent
	where $ J:= \left\lbrace  \chi{'} : \left| K_Y + L_{\chi{'}}\right| \ne \emptyset \right\rbrace  $ and $ R_{\sigma} $ is the reduced divisor supported on $ f^{*}\left( D_{\sigma} \right)  $.
	\end{Proposition}
    \noindent
    For the explanation of this proposition, we refer the reader to \cite[\rm Page 3]{2018arXiv180711854G}, \cite[\rm Remark 3.16]{MendesLopes2023} or \cite[\rm Section 3.4]{MR2067044}.
   
    \subsection{Resolution of singularities}
    Let $ Y $ be a smooth surface, and let $ \xymatrix{f: X \ar[r]& Y} $ be a smooth $  \mathbb{Z}_{2}^3 $-cover with building data $ \left\lbrace L_{\chi}, D_{\sigma} \right\rbrace_{\chi,\sigma}$. We consider $ L'_{\chi} \equiv  L_{\chi}$, $ D'_{\sigma} \equiv D_{\sigma} $ so that $ \left\lbrace L'_{\chi}, D'_{\sigma} \right\rbrace_{\chi,\sigma}$ is a building data of a $  \mathbb{Z}_{2}^3 $-cover $ \xymatrix{f: X' \ar[r]& Y} $ such that $ X' $ is a normal surface. Assume the moment that there is only one point $ P $ of $ B':= \sum\limits_{\sigma \ne 0} {D'_{\sigma} } $ giving rise to sigularities on $ X' $ and that the branch components $ D'_{\sigma} $ meet tranversally at $ P $. This section is devoted to presenting a resolution of these singularities, which is similar to the resolution of singular bidouble cover \cite{MR1718139}. The results in this section are taken from the work of C. Liedtke \cite[\rm Section 3]{MR2067044}.\\
    
   Let $\xymatrix{\pi: \tilde{Y} \ar[r] & Y}$ be the blow-up at $P$, and let $E$ be the exceptional divisor. We write $P = (\mu_{001}, \mu_{010}, \mu_{011}, \mu_{100}, \mu_{101}, \mu_{110}, \mu_{111})$, where $\mu_{\sigma}$ is the multiplicity of $D'_{\sigma}$ at $P$. Let $\tilde{\mu}_{\sigma}:= \mu_{\sigma}$ for all $\sigma$, except for at most one case where $\tilde{\mu}_{\sigma}:= \mu_{\sigma} -1$ so that the following equations has a solution $\left(a_{001}, a_{010}, a_{011}, a_{100}, a_{101}, a_{110}, a_{111}\right)$:
   $$
   \begin{adjustbox}{max width=\textwidth}
   	\begin{tabular}{l l l l l l l l l}
   		&&$  $&$ $&$ \tilde{\mu}_{100} $&$ +\tilde{\mu}_{101 } $&$ +\tilde{\mu}_{110} $&$ +\tilde{\mu}_{111} $& $ = 2a_{100}$\\
   		&&$ \tilde{\mu}_{010} $&$ +\tilde{\mu}_{011} $&$  $&$ $&$ +\tilde{\mu}_{110} $&$ +\tilde{\mu}_{111} $& $ = 2a_{010}$\\
   		&$ \tilde{\mu}_{001 } $&$ $&$ +\tilde{\mu}_{011} $&$$&$ +\tilde{\mu}_{101} $&$  $&$ +\tilde{\mu}_{111 } $& $ = 2a_{001}$\\
   		&&$ \tilde{\mu}_{010 } $&$ +\tilde{\mu}_{011} $&$ +\tilde{\mu}_{100} $&$ +\tilde{\mu}_{101} $&$  $&$  $&$ = 2a_{110}$\\
   		&$  \tilde{\mu}_{001} $&$ $&$ +\tilde{\mu}_{011} $&$ +\tilde{\mu}_{100} $&$  $&$ +\tilde{\mu}_{110} $&$  $& $ = 2a_{101}$\\
   		&$  \tilde{\mu}_{001} $&$ +\tilde{\mu}_{010} $&$  $&$  $&$ +\tilde{\mu}_{101} $&$ +\tilde{\mu}_{110} $&$  $& $ = 2a_{011}$\\
   		&$  \tilde{\mu}_{001} $&$ +\tilde{\mu}_{010} $&$  $&$ +\tilde{\mu}_{100} $&$ $&$ $&$ +\tilde{\mu}_{111 } $&$ =  2a_{111}$.
   	\end{tabular}
   \end{adjustbox}
   $$	
    
    \noindent
    We consider 
    \begin{align*}
    	\tilde{D}_{\sigma}:= &\pi^{*}\left( D'_{\sigma} \right) - \tilde{\mu}_{\sigma}E,\\
    	\tilde{L}_{\chi}:= &\pi^{*}\left( L'_{\chi} \right) - a_{\chi}E,
    \end{align*}
    
    \noindent
    for al non-trivial $ \sigma, \chi $. We notice that in the case where $ \tilde{\mu}_{\sigma}:= \mu_{\sigma} -1 $, we set $ \tilde{D}_{\sigma}:= \left( \pi^{*}\left( D'_{\sigma} \right) -\mu_{\sigma}E\right) + E $. This means that the exceptional divisor $ E $ is added into the strict transform of $ D'_{\sigma} $. Then $ \left\lbrace \tilde{L}_{\chi}, \tilde{D}_{\sigma} \right\rbrace_{\chi,\sigma}$ determines a canonical $  \mathbb{Z}_{2}^3 $-cover $ \xymatrix{\tilde{f}: \tilde{X} \ar[r]& \tilde{Y}} $. In \cite[\rm Section 3.2]{MR2067044}, the effect of point $P$ on $\chi(\mathcal{O}_X)$ and $K_X^2$ is thoroughly analyzed. To ensure a comprehensive understanding of our construction, we provide restatements of some of the cases we utilize. Specifically, our construction involves the utilization of two types of triple points, which we describe below.
    \begin{Proposition}\label{effect_triple_point}
    	Let $ \mu_{\sigma_1} = \mu_{\sigma_2} = \mu_{\sigma_3} = 1 $ and $ \mu_{\sigma} = 0 $ for all remaining $ \sigma $.
    	\begin{enumerate}
    		\item If $ \left\langle \sigma_1, \sigma_2, \sigma_3 \right\rangle \cong \mathbb{Z}_2^2 $, then $ K^2_{\tilde{X}}= K^2_{X} - 2$ and $ \chi(\mathcal{O}_{\tilde{X}}) = \chi(\mathcal{O}_X) $.
    		\item If $ \left\langle \sigma_1, \sigma_2, \sigma_3 \right\rangle \cong \mathbb{Z}_2^3 $, then $ K^2_{\tilde{X}}= K^2_{X} $ and $ \chi(\mathcal{O}_{\tilde{X}}) = \chi(\mathcal{O}_X) $.
    	\end{enumerate}
    \end{Proposition}
\begin{proof}
	1. By reordering indices, without loss of generality, we may assume that $ \sigma_1 = \left( 0,0,1\right), \sigma_2 = \left( 0,1,0\right), \sigma_3 = \left( 0,1,1\right) $. We have that
	\begin{align*}
		(\tilde{\mu}_{001}, \tilde{\mu}_{010}, \tilde{\mu}_{011}, \tilde{\mu}_{100}, \tilde{\mu}_{101}, \tilde{\mu}_{110}, \tilde{\mu}_{111})= \left( 1,1,1,0,0,0,0\right)\\
		(a_{001}, a_{010}, a_{011}, a_{100}, a_{101}, a_{110}, a_{111})= \left( 1,1,1,0,1,1,1\right).
	\end{align*}
	\noindent
	By Proposition \ref{invariants of Z_2^n cover} or \cite[\rm Section 3.2]{MR2067044}, we get $ K^2_{\tilde{X}}= K^2_{X} - 2$ and $ \chi(\mathcal{O}_{\tilde{X}}) = \chi(\mathcal{O}_X) $.\\
	\noindent
	2. Analogously, we may assume that $ \sigma_1 = \left( 0,0,1\right), \sigma_2 = \left( 0,1,0\right), \sigma_3 = \left( 1,0,0\right) $. We have that
	\begin{align*}
		(\tilde{\mu}_{001}, \tilde{\mu}_{010}, \tilde{\mu}_{011}, \tilde{\mu}_{100}, \tilde{\mu}_{101}, \tilde{\mu}_{110}, \tilde{\mu}_{111})= \left( 1,1,0,1,0,0,-1\right)\\
		(a_{001}, a_{010}, a_{011}, a_{100}, a_{101}, a_{110}, a_{111})= \left( 0,0,1,0,1,1,1\right).
	\end{align*}
	\noindent
	By Proposition \ref{invariants of Z_2^n cover} or \cite[\rm Section 3.2]{MR2067044}, we get $ K^2_{\tilde{X}}= K^2_{X}$ and $ \chi(\mathcal{O}_{\tilde{X}}) = \chi(\mathcal{O}_X) $.
	\end{proof}

	In our construction, we apply a type of quadruple points in order to decrease $ \chi(\mathcal{O}_X) $. These quadrupe points belong to the branch locus of a double cover ramified on $ 2L_{\chi} $ for some $ \chi $.
	\begin{Proposition}\label{Effect_quadrupe_point}
		Fix a character $ \chi $ and let $ \sigma_1, \sigma_2, \sigma_3, \sigma_4 $ be non-trivial elements in $ \mathbb{Z}_2^3 $ such that $ \chi\left( \sigma_1\right) = \chi\left( \sigma_2\right)= \chi\left( \sigma_3\right)= \chi\left( \sigma_4\right) = -1 $. Let $ \mu_{\sigma_1} = 2 $, $ \mu_{\sigma_2} = \mu_{\sigma_3} =1, \mu_{\sigma_4} = 0 $ and $ \mu_{\sigma} = 0 $ for all remaining $ \sigma $. Then 
		\begin{center}
			$ K^2_{\tilde{X}}= K^2_{X} - 2$ and $ \chi(\mathcal{O}_{\tilde{X}}) = \chi(\mathcal{O}_X) -1 $.
		\end{center}		
	\end{Proposition}
	\begin{proof}
	By reordering indices, without loss of generality, we may assume that $ \chi = \chi_{100} $ and $ \sigma_1 = \left( 1,0,0\right), \sigma_2 = \left( 1,0,1\right), \sigma_3 = \left( 1,1,0\right), \sigma_4 = \left( 1,1,1\right) $. We have that
	\begin{align*}
		(\tilde{\mu}_{001}, \tilde{\mu}_{010}, \tilde{\mu}_{011}, \tilde{\mu}_{100}, \tilde{\mu}_{101}, \tilde{\mu}_{110}, \tilde{\mu}_{111})= \left( 0,0,-1,2,1,1,0\right)\\
		(a_{001}, a_{010}, a_{011}, a_{100}, a_{101}, a_{110}, a_{111})= \left( 0,0,1,2,1,1,1\right).
	\end{align*}
	\noindent
	By Proposition \ref{invariants of Z_2^n cover} or \cite[\rm Section 3.2]{MR2067044}, we obtain $ K^2_{\tilde{X}}= K^2_{X} - 2$ and $ \chi(\mathcal{O}_{\tilde{X}}) = \chi(\mathcal{O}_X) -1$.
	\end{proof}

   	\section{Constructions of the regular surfaces}
   		In this section, we present the constructions of the regular surfaces with $ d = 10, 12, 14 $ listed in Theorem \ref{the theorem with d = 10 12 14}. We first construct a surface with $ d = 14 $ and $ q = 0 $ by using a $ \mathbb{Z}_2^3 $-cover of a del Pezzo surface of degree $ 7 $. To obtain regular surfaces with $ d = 10, 12 $, we impose triple points in the total branch locus and resolve the singularities.		   
   		\subsection{A surface with $ d = 14 $ and $ q = 0 $.}\label{the surface with d = 14 and q = 0} 
   		Let $ Y_2 $ be a del Pezzo surface of degree $ 7 $ (see Notation \ref{Notation of del Pezzo surface of degree 7}). We consider the following smooth divisors:
   		\begin{align*}
   		D_{010}:=& f_{11} + f_{12}  & D_{011}:=& f_{21} +e_1 \\
   		D_{100}:=&f_{22}            & D_{101}:=& C_1 \in \left| l +f_1\right| \\
   		D_{110}:=&C_2 \in \left| l+ f_1\right| & D_{111}:=&f_{23},
   		\end{align*}
   		\noindent
   		and $ D_{001}=0 $ where $ f_{11}, f_{12} \in \left| f_1\right| $, $ f_{21}, f_{22}, f_{23} \in \left| f_2\right| $ and $ C_1, C_2 $ are distinct divisors of $ Y_2 $ such that no more than two of these divisors $ D_{\sigma} $ go through the same point. We consider the following non-trivial divisors of $ Y_2 $:
   		$$	
   		\begin{tabular}{l r r r }
   		$ L_{100} :=$ &  $ f_1 $& $ +f_2 $&$ +l $\\
   		$ L_{010} :=$ &  $ f_1 $& $ +f_2 $&$ +l $\\
   		$ L_{001} :=$ &  $ $& $ f_2 $&$ +l$\\
   		$ L_{110} :=$ &  $ f_1 $& $ +f_2 $&$ +l $\\
   		$ L_{101} :=$ &  $ $& $ f_2 $&$ +l $\\
   		$ L_{011} :=$ &  $ 2f_1 $& $ $&$ +l $\\
   		$ L_{111} :=$ &  $ f_1 $& $ +f_2 $.&$ $
   		\end{tabular} 
   		$$
   		
   		\noindent
   		These divisors $ D_{\sigma}, L_{\chi} $ satisfy the following relations:
   		$$
   		\begin{adjustbox}{max width=\textwidth}
   		\begin{tabular}{l l l l l l l l l r r r}
   		$ 2L_{100} $&$ \equiv $&$  $&$ $&$ D_{100} $&$ +D_{101 } $&$ +D_{110} $&$ +D_{111} $ &$ \equiv $&  $ 2f_1 $& $ +2f_2 $&$ +2l $\\
   		$ 2L_{010} $&$ \equiv $&$ D_{010} $&$ +D_{011} $&$  $&$ $&$ +D_{110} $&$ +D_{111} $ &$ \equiv $&  $ 2f_1 $& $ +2f_2 $&$ +2l $\\
   		$ 2L_{001} $&$ \equiv D_{001 } $&$ $&$ +D_{011} $&$$&$ +D_{101} $&$  $&$ +D_{111 } $ &$ \equiv $&  $ $& $ 2f_2 $&$ +2l$ \\
   		$ 2L_{110} $&$ \equiv $&$ D_{010 } $&$ +D_{011} $&$ +D_{100} $&$ +D_{101} $&$  $&$  $& $ \equiv $&  $ 2f_1 $& $ +2f_2 $&$ +2l $\\
   		$ 2L_{101} $&$ \equiv D_{001} $&$ $&$ +D_{011} $&$ +D_{100} $&$  $&$ +D_{110} $&$  $ &$ \equiv $&  $ $& $ 2f_2 $&$ +2l $\\
   		$ 2L_{011} $&$ \equiv D_{001} $&$ +D_{010} $&$  $&$  $&$ +D_{101} $&$ +D_{110} $&$  $ &$ \equiv $&  $ 4f_1 $& $ $&$ +2l $\\
   		$ 2L_{111} $&$ \equiv D_{001} $&$ +D_{010} $&$  $&$ +D_{100} $&$ $&$ $&$ +D_{111 } $&$ \equiv $&  $ 2f_1 $& $ +2f_2 $.&
   	    \end{tabular}
        \end{adjustbox}
         $$	
   		
   		\noindent
   		Thus by Proposition \ref{Construction of cover of degree 8}, the divisors $ D_{\sigma}, L_{\chi} $ define a $ \mathbb{Z}^3_2 $-cover $ \xymatrix{g: X \ar[r] & Y_2}  $. Because each branch component $ D_{\sigma} $ is smooth and the total branch locus $ B $ is a normal crossings divisor, the surface $ X $ is smooth. Moreover, by Proposition \ref{invariants of Z_2^n cover}, the surface $ X $ satisfies the following:		
   		\begin{align*}
   		2K_X &\equiv g^*\left( f_1 + f_2 +l\right).
   		\end{align*}
   		\noindent
   		We notice that a surface is of general type and minimal if the canonical divisor is big and nef (see e.g. \cite[\rm Section 2]{MR2931875}). Since the divisor $ 2K_X $ is the pull-back of a nef and big divisor, the canonical divisor $ K_X $ is nef and big. Thus, the surface $ X $ is of general type and minimal. Furthermore, from Proposition \ref{invariants of Z_2^n cover}, the surface $X$ possesses the following invariants:
   		\begin{align*}
   		K_X^2=& 2\left( f_1 + f_2 +l\right)^2=14,\\
   		p_g\left( X\right) =&p_g\left( Y_2 \right) +\sum\limits_{\chi \ne  \chi_{000}  }{h^0\left( L_{\chi} + K_{Y_2} \right)}=h^0\left( \mathcal{O}_{Y_2} \right)+h^0\left( \mathcal{O}_{Y_2} \right)+h^0\left( \mathcal{O}_{Y_2} \right) =3,\\
   		\chi\left( \mathcal{O}_X\right) =&8\chi\left( \mathcal{O}_{Y_2} \right)  +\sum\limits_{\chi \ne \chi_{000}  }{\frac{1}{2}L_{\chi}\left( L_{\chi}+K_{Y_2}\right)}=4,\\
   		q\left( X\right)  =& 1+p_g\left( X\right)-\chi\left( \mathcal{O}_X\right)=0.   
   		\end{align*}
   		
   		We show that the canonical map $ \varphi_{\left| K_X  \right|}  $ has degree $ 14 $. By Proposition \ref{The generators of canonical system of Z_2^n cover}, the linear system $ \left| K_X\right|  $ is generated by the three following divisors:
   		\begin{align*}
   		\overline{f}_{11}+\overline{f}_{12} +\overline{f}_{21} + \overline{e}_{1}, \overline{f}_{22}+\overline{C}_{1}, \overline{C}_2+\overline{f}_{23},
   		\end{align*}
   		\noindent
   		where $ \overline{C}_1:=g^{*}\left( C_1\right)_{\text{red}}   $, $ \overline{C}_2:=g^{*}\left( C_2\right)_{\text{red}}   $, $ \overline{e}_1:=g^{*}\left( e_1\right)_{\text{red}}   $ and $ \overline{f}_{ij}:=g^{*}\left( f_{ij}\right)_{\text{red}}   $. Because the three divisors $ \overline{f}_{11}+\overline{f}_{12} +\overline{f}_{21} + \overline{e}_{1} $, $ \overline{f}_{22}+\overline{C}_{1} $, $ \overline{C}_2+\overline{f}_{23} $ have no common intersection, the linear system $ \left| K_X\right|  $ is base point free. This together with $ K_X^2 = 14 > 0 $ implies that the linear system $ \left| K_X\right|  $ is not composed with a pencil. Thus, the canonical image is $ \mathbb{P}^2 $ and the canonical map is of degree $ 14 $. Therefore, we obtain the surface described in the seventh row of Theorem \ref{the theorem with d = 10 12 14}.  
   		
   		\begin{Remark}
   		The surface $ X $ has two pencils of genus $ 5 $ corresponding to the fibres $ f_1, f_2 $.
   		\end{Remark}
   		\noindent
   		Infact, the $ \mathbb{Z}^3_2 $-cover $ \xymatrix{g: X \ar[r] & Y_2}  $ can be written as the composition of the three double covers: 
   		       $$ \xymatrix{X \ar[r]^{g_3} &X_2 \ar[r]^{g_2} & X_1 \ar[r]^{g_1} & Y_2}, $$
   		\noindent
   		where the double covers $g_1$, $g_2$, and $g_3$ correspond to $2L_{001}$, $2L_{010}$, and $2L_{100}$, respectively. The pull-back $ g^{*}\left( f_1\right) = g_3^{*}\left(  g_2^{*}\left(  g_1^{*}\left(  f_1 \right) \right) \right)  $ of a general fiber $ f_1 $ is an irreducible curve since the pull-backs of a general fiber $f_1$ intersect the branch locus in every double cover. By applying the Hurwitz formula, we find that the genus of $g^{*}(f_1)$ is $5$.

   		\subsection{Variations} 
   		We now consider variations of the construction presented in Section \ref{the surface with d = 14 and q = 0} in order to obtain the remaining regular surfaces listed in Theorem \ref{the theorem with d = 10 12 14}. To achieve this, we impose ordinary triple points on the branch locus and we resolve the singularities. We apply triple points that arise from the intersection of three divisors $D_{\sigma_1}, D_{\sigma_2},$ and $D_{\sigma_1 + \sigma_2}$. By introducing these triple points, we can achieve surfaces with the same Euler characteristic. According to Proposition \ref{effect_triple_point}, imposing one such triple point will result in a surface with $K^2 = 12$ and $\chi = 4$. On the other hand, if we simultaneously impose two such triple points, we will obtain a surface with $K^2 = 10$ and $\chi = 4$. More precisely, the specific types of triple points we impose and the corresponding surfaces obtained are described in the following table:
   			$$
   		\begin{tabular}{|c|c|c|c|c|c|}
   			\hline
   			$P = (\mu_{001}, \mu_{010}, \mu_{011}, \mu_{100}, \mu_{101}, \mu_{110}, \mu_{111})$& $ K_X^2  $&$ \chi\left( \mathcal{O}_X\right) $ &$ p_g\left( X\right) $ &$ q\left( X\right) $& $ d $ \\
   			\hline
   			$P_3 = (0, 0, 1, 0, 1, 1, 0)$&$ 12 $&$ 4 $&$ 3 $&$ 0 $& $ 12 $\\
   			\hline
   			$P_3 = (0, 0, 1, 0, 1, 1, 0)$ and $P_4 = (0, 1, 0, 0, 1, 0, 1)$&$ 10 $&$ 4 $&$ 3 $&$ 0 $& $ 10 $\\
   			\hline   			
   		\end{tabular} 
   		$$
   		\noindent
   		In our construction, we consistently have $ L_{100} \equiv L_{010} \equiv L_{110} \equiv -K_{Y_i} $, where $ Y_i $ is the del Pezzo we take the cover. This condition ensures that $ p_g\left( X\right)  = 3$. The degree of the canonical map is proven similarly in Section \ref{the surface with d = 14 and q = 0}. For further information regarding these constructions, we provide the building data of the $\mathbb{Z}_2^3$-covers.
   		
   		\subsubsection{The building data of the surface with $ d = 12 $ and $ q = 0 $}
   		Let us denote by $ Y_3 $ the blow-up of $ Y_2 $ at a point $ P_3 $. The surface $ Y_3 $ is a del Pezzo surface of degree $ 6 $ (see Notation \ref{Notation of del Pezzo surface of degree 6}). We consider the following smooth divisors:
   		\begin{align*}
   			D_{010}:=& f_{11} + f_{12}  & D_{011}:=& h_{23} +e_1 \\
   			D_{100}:=&f_{21}            & D_{101}:=& C_1 \in \left| f_1 + f_{3}\right| \\
   			D_{110}:=&C_2 \in \left| f_1+ f_3\right| & D_{111}:=&f_{22}
   		\end{align*}
   		\noindent
   		and $ D_{001}= 0 $, where $ f_{11}, f_{12} \in \left| f_1\right| $, $ f_{21}, f_{22} \in \left| f_2\right| $ and $ C_{1}, C_{2} $ are distinct divisors in $ Y_3 $ such that no more than two of these divisors $ D_{\sigma} $ go through the same point. In addition, we consider the following non-trivial divisors in $ Y_3 $:
   		$$	
   		\begin{tabular}{l r r r r}
   			$ L_{100} :=$ &  $ f_1 $& $ +f_2 $&$ +f_3 $&\\
   			$ L_{010} :=$ &  $ f_1 $& $ +f_2 $&$ +f_3 $&\\
   			$ L_{001} :=$ &  $  $& $ f_2 $&$ +f_3 $&\\
   			$ L_{110} :=$ &  $ f_1 $& $ +f_2 $&$ +f_3 $&\\
   			$ L_{101} :=$ &  $  $& $ f_2 $&$ +f_3 $&\\
   			$ L_{011} :=$ &  $ 2f_1 $& $ $&$ +f_3 $&$  $\\
   			$ L_{111} :=$ &  $ f_1 $& $ +f_2. $&$ $&
   		\end{tabular} 
   		$$
   		\noindent
   		These divisors $ D_{\sigma}, L_{\chi} $ define a $ \mathbb{Z}^3_2 $-cover $ \xymatrix{g: X \ar[r] & Y_3}  $. Therefore we obtain the surface in the fourth row of Theorem \ref{the theorem with d = 10 12 14}.
   		
   		\subsubsection{The building data of the surface with $ d = 10 $ and $ q = 0 $}
   		We denote by $ Y_4 $ the blow-up of $ Y_2 $ at two distinct points $ P_3 $ and $ P_4 $. The surface $ Y_4 $ is a del Pezzo surface of degree $ 5 $ (see Notation \ref{Notation of del Pezzo surface of degree 5}). We consider the following smooth divisors:
   		\begin{align*}
   		D_{010}:=& f_{11} + h_{14}  & D_{011}:=& h_{23} +e_1 \\
   		D_{100}:=&f_{21}             & D_{101}:=& C_1 \in \left| f_1 + h_{34}\right|  \\
   		D_{110}:=&C_2 \in \left| f_1 + f_{3}\right| & D_{111}:=&h_{24}
   		\end{align*}
   		\noindent
   		and $ D_{001} = 0 $, where $ f_{11} \in \left| f_1\right| $, $ f_{21} \in \left| f_2\right| $ and $ C_1, C_2 $  are divisors of $ Y_4 $ such that no more than two of these divisors $ D_{\sigma} $ go through the same point. We consider the following divisors:
   		$$	
   		\begin{tabular}{l r r r r r }
   			$ L_{100} :=$ &  $ f_1 $& $ +f_2 $&$ +f_3 $&$ -e_4 $&\\
   			$ L_{010} :=$ &  $ f_1 $& $ +f_2 $&$ +f_3 $&$ -e_4 $&\\
   			$ L_{001} :=$ &  $  $& $ f_2 $&$ +f_3 $&$ -e_4 $&\\
   			$ L_{110} :=$ &  $ f_1 $& $ +f_2 $&$ +f_3 $&$ -e_4 $&\\
   			$ L_{101} :=$ &  $ $& $ f_2 $&$ +f_3 $& $  $&\\
   			$ L_{011} :=$ &  $ 2f_1 $& $ $&$ +f_3 $&$ -e_4 $&\\
   			$ L_{111} :=$ &  $ f_1 $& $ +f_2 $&$ $&$ -e_4. $&
   		\end{tabular} 
   		$$
   		\noindent
   		These divisors $ D_{\sigma}, L_{\chi} $ define a $ \mathbb{Z}^3_2 $-cover $ \xymatrix{g: X \ar[r] & Y_4}  $. Therefore we obtain the surface in the first row of Theorem \ref{the theorem with d = 10 12 14}.
   		
   		\begin{Remark}
   		The surface $ X $ has four pencils of genus $ 5 $ corresponding to the fibres $ f_1, f_2, f_3, f_4$.
   		\end{Remark}
   		
   		\begin{Remark}
   		{\rm Taking the  $ \mathbb{Z}_2^3 $-cover  of  $ Y_2 $ ramified on the above branch locus with two triple points $ P_3, P_4 $, we would obtain a surface with four singular points of type $ \frac{1}{4}\left( 1,1\right) $ (cf. \cite[ \rm Table 1, Section 3.3]{MR2956036}). The preimage of each $ P_i $ is two $ \left( -4\right)  $-curves. The surface $ X $ is the minimal resolution of this singular surface.}
   		\end{Remark}               
   		\noindent
   		In fact, we consider the following commutative diagram
   		$$
   		\xymatrix{\overline{X} \ar[0,2]^{\mathbb{Z}^3_2}_{\overline{g}} && Y_2\\
   		X\ar[-1,0]^{c} \ar[0,2]^{\mathbb{Z}^3_2}_{g} && Y_4 \ar[-1,0]^{\pi}}
   		$$
        \noindent
        where the map $ \overline{g} $ is the  $ \mathbb{Z}_2^3 $-cover  of  $ Y_2 $ ramified on the above branch locus with two triple points $ P_3, P_4 $ and the map $ \pi $ is the blow-up at $ P_3 $, $ P_4 $. Each  divisor $ g^{*}\left( \pi\left( P_i\right) \right)  $ are two $ \left( -4\right)  $-curves. The map $ c $ is the contraction of these four $ \left( -4\right)  $-curves.

        \section{Constructions of the irregular surfaces}
        In this section, we present the constructions of the surfaces with $ d = 10,12,14 $ and $ q = 1,2 $ listed in Theorem \ref{the theorem with d = 10 12 14}. We first construct a surface with $ d = 14 $ and $ q = 1 $ by using a $ \mathbb{Z}_2^3 $-cover of a del Pezzo surface of degree $ 7 $. To obtain the other irregular surfaces described in Theorem \ref{the theorem with d = 10 12 14}, we impose triple or quadruple points in the total branch locus and resolve the singularities.
        
        \subsection{A surface with $ d = 14 $ and $ q = 1 $}\label{the surface with d = 14 and q = 1}
        The construction of the surface with $d=14$ and $q=1$ is essentially similar to the construction of the regular surface with $d=14$ described in Section \ref{the surface with d = 14 and q = 0}. We construct this surface as a $ \mathbb{Z}_2^3 $-cover of a del Pezzo surface $ Y_2 $ of degree $ 7 $. However, in order to achive $ q =1 $, we take the branch locus in such a way there is a quotient $ X/{\mathbb{Z}^2_2} $, which is a double cover of $ Y_2 $ ramified in four fibers $ f_2 $. More precisely, we consider the following smooth divisors:
        \begin{align*}
        D_{010}:=& f_{11} + f_{12}  & D_{011}:=& f_{21} +e_1 \\
        D_{100}:=&l_{1} + e_2            & D_{101}:=& f_{22} + h_{12} \\
        D_{110}:=&C_2 \in \left| l+ f_1\right| & D_{111}:=&f_{23},
        \end{align*}
        \noindent
        and $ D_{001} = 0 $, where $ f_{11}, f_{12} \in \left| f_1\right| $, $ f_{21}, f_{22}, f_{23} \in \left| f_2\right| $, $ l_1 \in \left| l\right| $ and $ C_2 $ are distinct divisors in $ Y_2 $ such that no more than two of these divisors $ D_{\sigma} $ go through the same point. We consider the following non-trivial divisors in $ Y_2 $:
        $$	
        \begin{tabular}{l r r r }
        $ L_{100} :=$ &  $ f_1 $& $ +f_2 $&$ +l $\\
        $ L_{010} :=$ &  $ f_1 $& $ +f_2 $&$ +l $\\
        $ L_{001} :=$ &  $ $& $ 2f_2 $&$ $\\
        $ L_{110} :=$ &  $ f_1 $& $ +f_2 $&$ +l $\\
        $ L_{101} :=$ &  $ $& $  $&$ 2l $\\
        $ L_{011} :=$ &  $ 2f_1 $& $ +f_2 $ $ $&$  $\\
        $ L_{111} :=$ &  $ f_1 $& $  $&$+l $.
        \end{tabular} 
        $$
        \noindent
        Since the divisors $ D_{\sigma}, L_{\chi} $ satisfy the relations in Proposition \ref{Construction of cover of degree 8}, these divisors $ D_{\sigma}, L_{\chi} $ define a $ \mathbb{Z}^3_2 $-cover $ \xymatrix{g: X \ar[r] & Y_2}  $. Because each branch component $ D_\sigma $ is smooth and the total branch locus is a normal crossings divisor, the surface $ X $ is smooth. Moreover, we have $q(X)=1$ based on the given building data of $2L_{001}=4f_2$, and $p_g(X)=3$ since $L_{100} \equiv L_{010} \equiv L_{110} \equiv -K_{Y_2}$. Infact, according to Proposition \ref{invariants of Z_2^n cover}, the surface $X$ possesses the following invariants:
        \begin{align*}
        	K_X^2= 14, p_g\left( X\right) = 3, \chi\left( \mathcal{O}_X\right) =3, q\left( X\right)  = 1.   
        \end{align*}
        
        \noindent
        Moreover, since 
        \begin{align*}
        2K_X &\equiv g^*\left( f_1 + f_2 +l\right),
        \end{align*}
        \noindent
        the canonical divisor $ K_X $ is nef and big. Thus, the surface $ X $ is of general type and minimal. 
        
        We prove that the canonical map $ \varphi_{\left| K_X  \right|}  $ is of degree $ 14 $. By Proposition \ref{The generators of canonical system of Z_2^n cover}, the linear system $ \left| K_X\right|  $ is generated by the three following divisors:
        \begin{align*}
        \overline{f}_{11}+\overline{f}_{12} +\overline{f}_{21} + \overline{e}_{1},     \overline{l}_{1}+\overline{e}_{2} + \overline{f}_{22} + \overline{h}_{12}, 
        \overline{C}_2+\overline{f}_{23},
        \end{align*}
        \noindent
        where $ \overline{l}_1:=g^{*}\left( l_1\right)_{\text{red}}  $, $ \overline{C}_2:=g^{*}\left( C_2\right)_{\text{red}}  $, $ \overline{e}_i:=g^{*}\left( e_i\right)_{\text{red}}  $, $ \overline{h}_{12}:=g^{*}\left( f_{12}\right)_{\text{red}}  $ and $ \overline{f}_{ij}:=g^{*}\left( f_{ij}\right)_{\text{red}}  $. Because the three divisors $ \overline{f}_{11}+\overline{f}_{12} +\overline{f}_{21} + \overline{e}_{1} $, $ \overline{l}_{1}+\overline{e}_{2} + \overline{f}_{22} + \overline{h}_{12} $, $ \overline{C}_2+\overline{f}_{23} $ have no common intersection, the linear system $ \left| K_X\right|  $ is base point free. This together with $ K_X^2 = 14 > 0 $ implies that the linear system $ \left| K_X\right|  $ is not composed with a pencil. Thus, the canonical image is $ \mathbb{P}^2 $ and the canonical map is of degree $ 14 $. Therefore we obtain the surface in the last row of Theorem \ref{the theorem with d = 10 12 14}.  
        
        \begin{Remark}
        The surface $ X $ has one pencil of genus $ 5 $ corresponding to the fibres $ f_1 $. The Albanese pencil of $ X $, which comes from the fibration $ \left| f_2\right|  $, has genus $ 3 $.
        \end{Remark}
    	\noindent
    	In fact, we consider the following commutative diagram
    	$$
    	\xymatrix{X \ar[0,2]^{\mathbb{Z}_2^3}_{g} \ar[1,1]_{\mathbb{Z}^2_2}^{g_1}& & Y_2\\
    		&Z \ar[-1,1]_{\mathbb{Z}_2}^{g_2}&}   	
    	$$
    	\noindent
    	where $ \xymatrix{g_2: Z \ar[r]& Y_2} $ is the double cover ramified on $2L_{001} \equiv 4f_2$. The intermediate surface $ Z $ has $ p_g\left( Z\right) =0 $ and $ q\left( Z\right) =1 $, indicating that the Albanese morphism of $Z$ maps its image to a curve. Consequently, applying the universal property, the image of the Albanese morphism of $X$ is also a curve. By the Hurwitz formula, we find that the genus of $ g^{*}\left( f_1 \right) $ is $ 5 $, while the genus of $ g^{*}\left( f_2 \right) $ is $ 3 $.
    
    	\subsection{Variations} 
    	We now consider variations of the construction presented in Section \ref{the surface with d = 14 and q = 1} in order to obtain the remaining irregular surfaces listed in Theorem \ref{the theorem with d = 10 12 14}. Initially, we utilize triple points that arise from the intersection of three divisors $D_{\sigma_1}, D_{\sigma_2},$ and $D_{\sigma_1 + \sigma_2}$. By introducing these triple points, we can get surfaces with the same Euler characteristic. According to Proposition \ref{effect_triple_point}, a such single triple point results in a surface with $K^2 = 12$ and $\chi = 3$. Conversely, the imposition of two such triple points simultaneously yields a surface with $K^2 = 10$ and $\chi = 3$. To decrease the Euler characteristic by one, we impose one quadrupe point as described in Proposition \ref{Effect_quadrupe_point}. By doing so, we can obtain a surface with $K^2 = 12$ and $\chi = 2$. Finally, we consider the application of triple points arising from the intersection of three divisors $D_{\sigma_1}, D_{\sigma_2},$ and $D_{\sigma_3}$ with $ \left\langle \sigma_1, \sigma_2, \sigma_3 \right\rangle \cong \mathbb{Z}_2^3 $. By Proposition \ref{effect_triple_point}, the imposition of a such single triple point results in a surface with the same Euler characteristic and self-intersection of canonical class. Moreover, by simultaneously imposing one such triple point and one quadruple point, we obtain a surface with $K^2 = 12$ and $\chi = 2$. The specific types of singular points we impose and the corresponding surfaces obtained are detailed in the following table:
    	$$
    	\begin{tabular}{|c|c|c|c|c|c|}
    		\hline
    		$P = (\mu_{001}, \mu_{010}, \mu_{011}, \mu_{100}, \mu_{101}, \mu_{110}, \mu_{111})$& $ K_X^2  $&$ \chi\left( \mathcal{O}_X\right) $ &$ p_g\left( X\right) $ &$ q\left( X\right) $& $ d $ \\
    		\hline
    		$P_3 = (0, 1, 0, 1, 0, 1, 0)$&$ 12 $&$ 3 $&$ 3 $&$ 1 $& $ 12 $\\
    		\hline
    		$P_3 = (0, 1, 0, 1, 0, 1, 0)$ and $P_4 = (0, 1, 0, 1, 0, 1, 0)$&$ 10 $&$ 3 $&$ 3 $&$ 1 $& $ 10 $\\
    		\hline 
    			$P_3 = (0, 0, 1, 1, 0, 2, 0)$&$ 12 $&$ 2 $&$ 3 $&$ 2 $& $ 12 $\\
    		\hline
    		$P_3 = (0, 0, 1, 1, 0, 2, 0)$ and $P_4 = (0, 1, 0, 1, 0, 0, 1)$&$ 12 $&$ 2 $&$ 3 $&$ 2 $& $ 10 $\\
    		\hline   	  			
    	\end{tabular} 
    	$$
    	\noindent
    	In our construction, we consistently maintain $ L_{100} \equiv L_{010} \equiv L_{110} \equiv -K_{Y_i} $, where $ Y_i $ is the del Pezzo we take the cover. This condition guarantees that $ p_g\left( X\right)  = 3$. With the exception of the surface with $ d = 10, q = 2 $, the degree of the canonical map can be proven in a similar manner as described in Section \ref{the surface with d = 14 and q = 1}. The exceptional case will prove in the last section. For further details regarding these constructions, we provide the building data of the $\mathbb{Z}_2^3$-covers.
        
        \subsubsection{The building data of the surface with $ d = 12 $ and $ q = 1 $}         
        Let us denote by $ Y_3 $ the blow-up of $ Y_2 $ at a point $ P_3 $. The surface $ Y_3 $ is a del Pezzo surface of degree $ 6 $ (see Notation \ref{Notation of del Pezzo surface of degree 6}). We consider the following smooth divisors:
        \begin{align*}
        D_{010}:=& f_{11} + h_{13}  & D_{011}:=& f_{21} +e_1 \\
        D_{100}:=&f_{31} + e_2            & D_{101}:=& f_{22} + h_{12} \\
        D_{110}:=&C_2 \in \left| f_1+ f_3\right| & D_{111}:=&f_{23}
        \end{align*}
        \noindent
        and $ D_{001}= 0 $, where $ f_{11} \in \left| f_1\right| $, $ f_{21}, f_{22}, f_{23} \in \left| f_2\right| $, $ f_{31} \in \left| f_3\right| $ and $ C_2 $ are distinct divisors in $ Y_3 $ such that no more than two of these divisors $ D_{\sigma} $ go through the same point. In addition, we consider the following non-trivial divisors in $ Y_3 $:
        $$	
        \begin{tabular}{l r r r r}
        	$ L_{100} :=$ &  $ f_1 $& $ +f_2 $&$ +f_3 $&\\
        	$ L_{010} :=$ &  $ f_1 $& $ +f_2 $&$ +f_3 $&\\
        	$ L_{001} :=$ &  $  $& $ 2f_2 $&$  $&\\
        	$ L_{110} :=$ &  $ f_1 $& $ +f_2 $&$ +f_3 $&\\
        	$ L_{101} :=$ &  $ l $& $  $&$ +f_3 $&\\
        	$ L_{011} :=$ &  $ 2f_1 $& $+f_2 $&$  $&$ -e_3 $\\
        	$ L_{111} :=$ &  $ f_1 $& $  $&$+f_3 $.&
        \end{tabular} 
        $$
        \noindent
        These divisors $ D_{\sigma}, L_{\chi} $ define a $ \mathbb{Z}^3_2 $-cover $ \xymatrix{g: X \ar[r] & Y_3}  $. Moreover, the surface $X$ possesses the following invariants:
        \begin{align*}
        	K_X^2= 12, p_g\left( X\right) = 3, \chi\left( \mathcal{O}_X\right) =3, q\left( X\right)  = 1, d = 12  
        \end{align*}
                
        \begin{Remark}
        The surface $ X $ has two pencils of genus $ 5 $ corresponding to the fibres $ f_1 $, $ f_3 $. The Albanese pencil of $ X $, which comes from the fibration $ \left| f_2\right|  $, has genus $ 3 $.
        \end{Remark}
        
        \begin{Remark}
        {\rm Taking the  $ \mathbb{Z}_2^3 $ cover of $ Y_2 $ ramified on the above branch locus with a triple point $ P_3 $, we would obtain a surface with two singular points of type $ \frac{1}{4}\left( 1,1\right) $ coming from the point $ P_3 $ (cf. \cite[ \rm Table 1, Section 3.3]{MR2956036}). The surface $ X $ is the minimal resolution of this singular surface.}
        \end{Remark}
        
        \subsubsection{The building data of the surface with $ d = 10 $ and $ q = 1 $}         
        Let us denote by $ Y_4 $ the blow-up of $ Y_2 $ at two distinct points $ P_3 $ and $ P_4 $. The surface $ Y_4 $ is a del Pezzo surface of degree $ 5 $ (see Notation \ref{Notation of del Pezzo surface of degree 5}). We consider the following smooth divisors of $ Y_4 $:
        \begin{align*}
        	D_{010}:=& h_{13} + h_{14}  & D_{011}:=& f_{21} +e_1 \\
        	D_{100}:=&h_{34}  +e_2           & D_{101}:=& h_{12} + f_{22}  \\
        	D_{110}:=&C_2 \in \left| f_1 + h_{34}\right| & D_{111}:=&f_{23},
        \end{align*}
        \noindent
        and $ D_{001} = 0 $, where  $ f_{21}, f_{22}, f_{23} \in \left| f_2\right| $ and $ C_2 $ are distinct divisors such that no more than two of these divisors $ D_{\sigma} $ go through the same point. We consider the following non-trivial divisors of $ Y_4 $:
        $$	
        \begin{tabular}{l r r r r}
        	$ L_{100} \equiv$ &  $ f_1 $& $ +f_2 $&$ +f_3 $&$ -e_4 $\\
        	$ L_{010} \equiv$ &  $ f_1 $& $ +f_2 $&$ +f_3 $&$ -e_4 $\\
        	$ L_{001} \equiv$ &  $  $& $ 2f_2 $&$  $&$  $\\
        	$ L_{110} \equiv$ &  $ f_1 $& $ +f_2 $&$ +f_3 $&$ -e_4 $\\
        	$ L_{101} \equiv$ &  $ $& $  $&$ f_3 $& $ +f_4 $\\
        	$ L_{011} \equiv$ &  $ 2f_1 $& $+f_2 $&$ -e_3 $&$ -e_4 $\\
        	$ L_{111} \equiv$ &  $ f_1 $& $  $&$+f_3 $&$ -e_4. $
        \end{tabular} 
        $$
        \noindent
        These divisors $ D_{\sigma}, L_{\chi} $ define a $ \mathbb{Z}^3_2 $-cover $ \xymatrix{g: X \ar[r] & Y_4}  $. Moreover, the surface $X$ possesses the following invariants:
        \begin{align*}
        	K_X^2= 10, p_g\left( X\right) = 3, \chi\left( \mathcal{O}_X\right) =3, q\left( X\right)  = 1, d =10.   
        \end{align*}
                
        \begin{Remark}
        	The surface $ X $ has three pencils of genus $ 5 $ corresponding to the fibre $ f_1, f_3, f_4$. The Albanese pencil of $ X $, which comes from the fibration $ \left| f_2\right|  $, has genus $ 3 $.
        \end{Remark}
        
        \begin{Remark}
        	{\rm Taking the  $ \mathbb{Z}_2^3 $-cover  of  $ Y_2 $ ramified on the above branch locus with two ordinary triple points $ P_3, P_4 $, we would obtain a surface with four points of type $ \frac{1}{4}\left( 1,1\right) $ coming from $ P_3 $ and $ P_4 $   (cf. \cite[ \rm Table 1, Section 3.3]{MR2956036}). The surface $ X $ is the minimal resolution of this singular surface.}
        \end{Remark}

        \subsubsection{The building data of the surface with $ d = 12 $ and $ q = 2 $} 
        We denote by $ Y_3 $ the blow-up of $Y_2 $ at a point $ P_3 $. The surface $ Y_3 $ is a del Pezzo surface of degree $ 6 $ (see Notation \ref{Notation of del Pezzo surface of degree 6}). We consider the following divisors of $ Y_3 $:
        \begin{align*}
        D_{010}:=& f_{11} + f_{12}  & D_{011}:=& h_{23} +e_1 \\
        D_{100}:=&f_{31} + e_2            & D_{101}:=& f_{21} + h_{12} \\
        D_{110}:=&f_{32} + h_{13} & D_{111}:=&f_{22}+e_3,
        \end{align*}
        \noindent
        and $ D_{001} = 0 $, where $ f_{11}, f_{12} \in \left| f_1\right| $, $ f_{21}, f_{22} \in \left| f_2\right| $ and $ f_{31}, f_{32} \in \left| f_3\right| $ are distinct divisors such that no more than two of these divisors $ D_{\sigma} $ go through the same point. We consider the following non-trivial divisors of $ Y_3 $:
        $$	
        \begin{tabular}{l r r r r }
        $ L_{100} :=$ &  $ f_1 $& $ +f_2 $&$ +f_3 $&\\
        $ L_{010} :=$ &  $ f_1 $& $ +f_2 $&$ +f_3 $&\\
        $ L_{001} :=$ &  $  $& $ 2f_2 $&$  $&$  $\\
        $ L_{110} :=$ &  $ f_1 $& $ +f_2 $&$ +f_3 $&\\        
        $ L_{101} :=$ &  $  $& $  $&$ 2f_3 $&\\
        $ L_{011} :=$ &  $ 2f_1 $& $+f_2 $&$  $&$ -e_3 $\\
        $ L_{111} :=$ &  $ f_1 $& $  $&$ $&$ +l $.
        \end{tabular} 
        $$
        \noindent
        These divisors $ D_{\sigma}, L_{\chi} $ define a $ \mathbb{Z}^3_2 $-cover $ \xymatrix{g: X \ar[r] & Y_3}  $. Moreover, the surface $X$ possesses the following invariants:
        \begin{align*}
        K_X^2= 12, p_g\left( X\right) = 3, \chi\left( \mathcal{O}_X\right) =2, q\left( X\right)  = 2, d = 12.   
        \end{align*}
        
        \begin{Remark}
        {\rm Taking the  $ \mathbb{Z}_2^3 $-cover of $ Y_2 $ ramified on the above branch locus with a quadruple point $ P_3 $, we would obtain a surface with a singular point coming from the point $ P_3 $. This point is Gorenstein elliptic singularity whose minimal resolution are elliptic curves with self-intersection $-2$ (cf. \cite[ \rm Table 1, Section 3.3]{MR2956036}). The surface $ X $ is the minimal resolution of this singular surface.}
        \end{Remark}           
        
        \begin{Remark}
        {\rm The Albanese morphism of this surface is of degree $ 2 $.}
        \end{Remark}
        \noindent
        In fact, we consider the following commutative diagram
        $$
        \xymatrix{X \ar[0,2]^{\mathbb{Z}_2^3}_{g} \ar[1,1]_{\mathbb{Z}_2}^{g_1}& & Y_3\\
        	&Z \ar[-1,1]_{\mathbb{Z}_2^2}^{g_2}&}   	
        $$
        \noindent
        where $ \xymatrix{g_2: Z \ar[r]& Y_3} $ is the bidouble cover with the following buiding data:
        $$
        \begin{tabular}{l l l}
        	$ D_1:= D_{011} $ &$= h_{23} + e_1 $ &$\equiv l + e_1 - e_2 - e_3 $\\
        	$ D_2:= D_{101} + D_{111} $ & $ = f_{21} + h_{12} + f_{22} + e_3 $& $ \equiv 3l - e_1 - 3e_2 + e_3 $\\
        	$ D_3:= D_{100} + D_{110} $ & $ = f_{31} + e_2 + f_{32} + h_{13} $&$ \equiv 3l - e_1 + e_2 - 3e_3 $\\
        	$ L_1:= L_{100}  $&&$ \equiv 3l - e_1 - e_2 - e_3 $\\
        	$ L_2:= L_{101} $&&$\equiv 2l-2e_3 $\\
        	$ L_3:= L_{001} $&&$\equiv 2l-2e_2 $.   		
        \end{tabular}
        $$
        \noindent
        We have that $ p_g\left( Z\right) = 1, q\left( Z\right) =2 $ and 
        \begin{align*}
        2K_Z \equiv g_2^{*}\left( h_{23 }+ e_1\right). 
        \end{align*}
        
        \noindent
        We notice that the surface $ Z $ contains $ \left( -1\right)  $-curves. Let $ \xymatrix{c: Z \ar[r]& \overline{Z}} $ be the contraction map. The minimal surface $ \overline{Z} $ satisfies $ p_g\left( \overline{Z}\right) = 1, q\left( \overline{Z}\right) =2 $ and $ 2K_{\overline{Z}} \equiv 0$. So the surface $ \overline{Z} $ is an abelian surface. By the universal property of the Albanese morphism of the surface $ X $, the following diagram commutes:
        $$
        \xymatrix{X \ar[r]^{\mathbb{Z}_2} \ar[1,1]^{\alpha}& Z \ar[r]^{c} & \overline{Z}\\
        	&Alb\left( X\right) \ar[-1,1]&}
        $$
        \noindent
        Because the surface $ X $ is of general type, the Albanese morphism $ \alpha $ is of degree $ 2 $ and $ Alb\left( X\right)  $ is isomorphic to $ \overline{Z} $.

        \subsubsection{The building data of the surface with $ d = 10 $ and $ q = 2 $.}         
        Let us denote by $ Y_4 $ the blow-up of $ Y_2 $ at two distinct points $ P_3 $ and $ P_4 $. $ Y_4 $ is a del Pezzo surface of degree $ 5 $ (see Notation \ref{Notation of del Pezzo surface of degree 5}). We consider the following smooth divisors of $ Y_4 $:
        \begin{align*}
        D_{001}:=& e_{4}  &&  \\
        D_{010}:=& h_{14} + f_{11}  & D_{011}:=& h_{23} +e_1 \\
        D_{100}:=&h_{34} + e_{2}            & D_{101}:=& h_{12} + f_{21}  \\
        D_{110}:=&h_{13} + f_{31} & D_{111}:=&h_{24}+e_3,
        \end{align*}
        \noindent
        where $ f_{11} \in \left| f_1\right| $, $ f_{21} \in \left| f_2\right| $, $ f_{31} \in \left| f_3\right| $ are divisors such that no more than two of these divisors $ D_{\sigma} $ go through the same point. We consider the following non-trivial divisors of $ Y_4 $:
        $$	
        \begin{tabular}{l r r r r }
        $ L_{100} :=$ &  $ f_1 $& $ +f_2 $&$ +f_3 $&$ -e_4 $\\
        $ L_{010} :=$ &  $ f_1 $& $ +f_2 $&$ +f_3 $&$ -e_4 $\\
        $ L_{001} :=$ &  $  $& $ 2f_2 $&$  $&$  $\\
        $ L_{110} :=$ &  $ f_1 $& $ +f_2 $&$ +f_3 $&$ -e_4 $\\
        $ L_{101} :=$ &  $ $& $  $&$ 2f_3 $& $  $\\
        $ L_{011} :=$ &  $ 2f_1 $& $  +f_2$&$-e_3 $&$$\\
        $ L_{111} :=$ &  $ f_1 $& $ $&$  $&$ +f_4. $
        \end{tabular} 
        $$
        \noindent
        These divisors $ D_{\sigma}, L_{\chi} $ define a $ \mathbb{Z}^3_2 $-cover $ \xymatrix{g: X \ar[r] & Y_4}  $. Moreover, the surface $X$ possesses the following invariants:
        	\begin{align*}
        	K_X^2= 12, p_g\left( X\right) = 3, \chi\left( \mathcal{O}_X\right) =2, q\left( X\right)  = 2.   
        	\end{align*}
        	
        	We show that the canonical map $ \varphi_{\left| K_X  \right|}  $ has degree $ 10 $ and the linear system $ \left| K_X\right|  $ has a non-trivial fixed component. By Proposition \ref{The generators of canonical system of Z_2^n cover}, the linear system $ \left| K_X\right|  $ is generated by the three following divisors:
        	\begin{align*}
        	\overline{e}_{4}+\overline{h}_{14}+\overline{f}_{11} +\overline{h}_{23} + \overline{e}_{1}, 
        	\overline{e}_{4}+\overline{h}_{34}+ \overline{e}_{2}+\overline{h}_{12}+\overline{f}_{21}, 
        	\overline{e}_{4}+\overline{h}_{13}+ \overline{f}_{31}+\overline{h}_{24}+\overline{e}_{3},
        	\end{align*}
        	\noindent
        	where $ \overline{C}_i:=g^{*}\left( C_i\right)_{\text{red}}   $, $ \overline{e}_i:=g^{*}\left( e_i\right)_{\text{red}}   $, $ \overline{h}_{ij}:=g^{*}\left( h_{ij}\right)_{\text{red}}   $ and $ \overline{f}_{ij}:=g^{*}\left( f_{ij}\right)_{\text{red}}   $. Because the $ \left( -2\right) $-curve $\overline{e}_{4} $ is the common part of these three above divisors, the $ \left( -2\right) $-curve $\overline{e}_{4} $ is a fixed component of $ \left| K_X\right|  $.
        	
        	On the other hand, since the three divisors $\overline{h}_{14}+\overline{f}_{11} +\overline{h}_{23} + \overline{e}_{1} $, $ \overline{h}_{34}+ \overline{e}_{2}+\overline{h}_{12}+\overline{f}_{21} $, $ \overline{h}_{13}+ \overline{f}_{31}+\overline{h}_{24}+\overline{e}_{3} $ have no common intersection, the linear system $ \left| M\right|  $ is base point free, where $ \left| M\right| := \left| \overline{h}_{14}+\overline{f}_{11} +\overline{h}_{23} + \overline{e}_{1}\right| $. This together with $ M^2 = 10 > 0 $ implies that the linear system $ \left| K_X\right|  $ is not composed with a pencil. Thus, the canonical image is $ \mathbb{P}^2 $ and the canonical map is of degree $ 10 $. Therefore we obtain the surface described in the third row of Theorem \ref{the theorem with d = 10 12 14}.
        	
        	\begin{Remark}
        		{\rm Taking the  $ \mathbb{Z}_2^3 $-cover of $ Y_1 $ ramified on the above branch locus with singular points $ P_3, P_4 $, we would obtain a surface with one singular point of type $ A_1 $ and one Gorenstein elliptic singular point whose minimal resolution is an elliptic curve with self-intersection $-2$. The $ A_1 $ point comes from $  P_4 $ and the Gorenstein elliptic singular point comes from $ P_3 $. The surface $ X $ is the minimal resolution of this singular surface.}
        	\end{Remark}

        	\begin{Remark}
        	{\rm Similarly to the construction of the surface with $ d= 2, p_g =3, q = 2$, the Albanese morphism of this surface is of degree $ 2 $.}
        	\end{Remark}   
        	
        	\begin{Remark}
        	Let $ H $ be the subgroup generated by $ \left( 0,0,1\right)  $ of the group $ \mathbb{Z}_2^3 $. In the all above constructions, the canonical map $ \varphi_{\left| K_X  \right|} $ is the composition of the quotient map $ \xymatrix{X \ar[r]&X/H} $ with the canonical map of the quotient surface $ X/H $ because $ h^{0}\left( L_{\chi} + K_{Y}\right) = 0  $ for all $ \chi \in \left\lbrace  \chi_{001}, \chi_{101}, \chi_{011}, \chi_{111}  \right\rbrace  $.
        	\end{Remark}

\section*{Acknowledgments}
The author is deeply indebted to Margarida Mendes Lopes for all her help. Furthermore, the author would like to extend his heartfelt appreciation to the referee for providing numerous insightful suggestions and detailed comments that significantly improved the clarity of the paper. The author was partially supported by Funda\c{c}\~{a}o para a Ci\^{e}ncia e Tecnologia (FCT), Portugal through the program Lisbon Mathematics PhD (LisMath) of the University of Lisbon, scholarship FCT - PD/BD/113632/2015 and project UID/MAT/04459/2019 of CAMGSD. This paper was finished during the author's postdoctoral fellowship at the National Center for Theoretical Sciences (NCTS), Taiwan, under the grant number MOST 110-2123-M-002-005. The author would like to thank NCTS for the kind hospitality.


\begin{thebibliography}{10}
	
	\bibitem{MR2956036}
	{\sc Alexeev, V., and Pardini, R.}
	\newblock Non-normal abelian covers.
	\newblock {\em Compos. Math. 148}, 4 (2012), 1051--1084.
	
	\bibitem{MR4278662}
	{\sc Bauer, I., and Pignatelli, R.}
	\newblock Rigid but not infinitesimally rigid compact complex manifolds.
	\newblock {\em Duke Math. J. 170}, 8 (2021), 1757--1780.
	
	\bibitem{MR553705}
	{\sc Beauville, A.}
	\newblock L'application canonique pour les surfaces de type g\'en\'eral.
	\newblock {\em Invent. Math. 55}, 2 (1979), 121--140.	
	
	\bibitem{MR1718139}
	{\sc Catanese, F.}
	\newblock Singular bidouble covers and the construction of interesting
	algebraic surfaces.
	\newblock In {\em Algebraic geometry: {H}irzebruch 70 ({W}arsaw, 1998)},
	vol.~241 of {\em Contemp. Math.} Amer. Math. Soc., Providence, RI, 1999,
	pp.~97--120.
	
	\bibitem{2022arXiv220906057F}
	{\sc {Fallucca}, F.}
	\newblock {Examples of surfaces with canonical maps of degree $12$, $13$, $15$,
		$16$ and $18$}.
	\newblock {\em arXiv e-prints\/} (Sept. 2022), arXiv:2209.06057.
	
	\bibitem{2022arXiv220702969F}
	{\sc {Fallucca}, F., and {Gleissner}, C.}
	\newblock {Some surfaces with canonical maps of degree $10$, $11$ and $14$}.
	\newblock {\em arXiv e-prints\/} (July 2022), arXiv:2207.02969.
	
	\bibitem{2018arXiv180711854G}
	{\sc {Gleissner}, C., {Pignatelli}, R., and {Rito}, C.}
	\newblock {New surfaces with canonical map of high degree}.
	\newblock {\em arXiv e-prints\/} (July 2018), arXiv:1807.11854.
	
	\bibitem{MR4298919}
	{\sc Lai, C.-J., and Yeung, S.-K.}
	\newblock Examples of surfaces with canonical map of maximal degree.
	\newblock {\em Taiwanese J. Math. 25}, 4 (2021), 699--716.
	
	\bibitem{MR2067044}
	{\sc Liedtke, C.}
	\newblock Singular abelian covers of algebraic surfaces.
	\newblock {\em Manuscripta Math. 112}, 3 (2003), 375--390.
	
	\bibitem{MR2931875}
	{\sc Lopes, M.~M., and Pardini, R.}
	\newblock The geography of irregular surfaces.
	\newblock In {\em Current developments in algebraic geometry}, vol.~59 of {\em
		Math. Sci. Res. Inst. Publ.} Cambridge Univ. Press, Cambridge, 2012,
	pp.~349--378.
	
	\bibitem{MendesLopes2023}
	{\sc Mendes~Lopes, M., and Pardini, R.}
	\newblock {\em On the Degree of the Canonical Map of a Surface of General
		Type}.
	\newblock Springer International Publishing, Cham, 2023, pp.~305--325.
	
	\bibitem{MR4008073}
	{\sc Nguyen Bin}
	\newblock A new example of an algebraic surface with canonical map of degree
	16.
	\newblock {\em Arch. Math. (Basel) 113}, 4 (2019), 385--390.
	
	\bibitem{MR4203293}
	{\sc Nguyen Bin}
	\newblock Some infinite sequences of canonical covers of degree 2.
	\newblock {\em Adv. Geom. 21}, 1 (2021), 143--148.
	\bibitem{MR4334862}
	{\sc Nguyen Bin}
	\newblock Some examples of algebraic surfaces with canonical map of degree 20.
	\newblock {\em C. R. Math. Acad. Sci. Paris 359\/} (2021), 1145--1153.
	
	\bibitem{MR4394091}
	{\sc Nguyen Bin}
	\newblock New examples of canonical covers of degree 3.
	\newblock {\em Math. Nachr. 295}, 3 (2022), 450--467.
	
	\bibitem{MR1103912}
	{\sc Pardini, R.}
	\newblock Abelian covers of algebraic varieties.
	\newblock {\em J. Reine Angew. Math. 417\/} (1991), 191--213.
	
	\bibitem{MR1103913}
	{\sc Pardini, R.}
	\newblock Canonical images of surfaces.
	\newblock {\em J. Reine Angew. Math. 417\/} (1991), 215--219.
	
	\bibitem{MR527234}
	{\sc Persson, U.}
	\newblock Double coverings and surfaces of general type.
	\newblock In {\em Algebraic geometry ({P}roc. {S}ympos., {U}niv. {T}roms\o ,
		{T}roms\o , 1977)}, vol.~687 of {\em Lecture Notes in Math.} Springer,
	Berlin, 1978, pp.~168--195.
	
	\bibitem{MR3391024}
	{\sc Rito, C.}
	\newblock New canonical triple covers of surfaces.
	\newblock {\em Proc. Amer. Math. Soc. 143}, 11 (2015), 4647--4653.
	
	\bibitem{MR3663791}
	{\sc Rito, C.}
	\newblock A surface with canonical map of degree 24.
	\newblock {\em Internat. J. Math. 28}, 6 (2017), 1750041, 10.
	
	\bibitem{MR3619737}
	{\sc Rito, C.}
	\newblock A surface with {$q=2$} and canonical map of degree 16.
	\newblock {\em Michigan Math. J. 66}, 1 (2017), 99--105.
	
	\bibitem{MR4372410}
	{\sc Rito, C.}
	\newblock Surfaces with canonical map of maximum degree.
	\newblock {\em J. Algebraic Geom. 31}, 1 (2022), 127--135.
	
	\bibitem{MR1141782}
	{\sc Tan, S.~L.}
	\newblock Surfaces whose canonical maps are of odd degrees.
	\newblock {\em Math. Ann. 292}, 1 (1992), 13--29.
	
	\bibitem{MR842626}
	{\sc Xiao, G.}
	\newblock Algebraic surfaces with high canonical degree.
	\newblock {\em Math. Ann. 274}, 3 (1986), 473--483.
	
	\bibitem{2015arXiv151006622Y}
	{\sc {Yeung}, S.-K.}
	\newblock {A surface of maximal canonical degree}.
	\newblock {\em arXiv e-prints\/} (Oct. 2015), arXiv:1510.06622.
	
\end{thebibliography}

\CurrentAddresses

\Addresses

\SecondAddresses

\end{document}